\documentclass[10pt, letterpaper]{amsart}
\usepackage{mathbbol}
\usepackage{amsmath, amsthm, amssymb}
\usepackage{txfonts, mathrsfs}
\usepackage[usenames]{color}

\usepackage{psfrag}
\usepackage{graphicx}

\newtheorem{thm}{Theorem}[section]
\newcommand{\bt}{\begin{thm}}
\newcommand{\et}{\end{thm}}

\newtheorem{cor}[thm]{Corollary}
\newcommand{\bc}{\begin{cor}}
\newcommand{\ec}{\end{cor}}

\newtheorem{lem}[thm]{Lemma}
\newcommand{\bl}{\begin{lem}}
\newcommand{\el}{\end{lem}}

\newtheorem{prop}[thm]{Proposition}
\newcommand{\bp}{\begin{prop}}
\newcommand{\ep}{\end{prop}}

\newtheorem{defn}[thm]{Definition}
\newcommand{\bd}{\begin{defn}}      
\newcommand{\ed}{\end{defn}}

\newtheorem{rmrk}[thm]{Remark}
\newcommand{\br}{\begin{rmrk}}
\newcommand{\er}{\end{rmrk}}

\newtheorem{quest}[thm]{Question}
\newcommand{\bq}{\begin{quest}}
\newcommand{\eq}{\end{quest}}

\newcommand{\thmref}[1]{Theorem~\ref{#1}}
\newcommand{\secref}[1]{Section~\ref{#1}}

\newcommand{\propref}[1]{Proposition~\ref{#1}}

\newcommand{\N}{\mathbb{N}}

\newcommand{\R}{\mathbb{R}}
\newcommand{\Z}{\mathbb{Z}}

\newcommand{\p}{\partial}
\newcommand{\spa}{\operatorname{span}}

\newcommand{\diam}{\operatorname{diam}}
\newcommand{\ANdim}{\operatorname{dim}_{AN}}

\begin{document}

\title{Lipschitz extensions into Jet space Carnot groups}

\author{Stefan Wenger}

\address
  {Department of Mathematics\\
University of Illinois at Chicago\\
851 S. Morgan Street\\
Chicago, IL 60607--7045 }
\email{wenger@math.uic.edu}

\author{Robert Young}

\address{Institut des Hautes \'Etudes Scientifiques\\
Le Bois Marie\\
35 route de Chartres\\
F--91440 Bures-sur-Yvette, France }
\email{rjyoung@ihes.fr}

\date{\today}

\thanks{The first author was partially supported by NSF grant DMS 0956374}

\begin{abstract}
 The aim of this article is to prove a Lipschitz extension theorem for partially defined Lipschitz maps to jet spaces endowed with a left-invariant sub-Riemannian Carnot-Carath\'eodory distance.
 The jet spaces give a model for a certain class of Carnot groups, including in particular all Heisenberg groups. The proofs or our theorems are rather elementary, they are inspired by ideas in papers of Gromov, Young, and Lang-Schlichenmaier. 
\end{abstract}

\maketitle

\bigskip

\section{Introduction}

Let $X$ and $Y$ be metric spaces. Recall that a map $f: X\to Y$ is said to be $\lambda$-Lipschitz if
\begin{equation*}
 d_Y(f(x), f(x'))\leq \lambda\, d_X(x,x')
\end{equation*}
for all $x,x'\in X$, where $d_X$ and $d_Y$ denote the metrics on $X$ and $Y$, respectively. A natural question asks: for which (classes of) spaces $X, Y$ does every Lipschitz map $f: Z\to Y$, defined on a subset $Z\subset X$, admit a Lipschitz extension $\bar{f}: X\to Y$ to all of $X$?
This and related problems have attracted research interest for many
years. The earliest results dealt with $Y=\R^n$, among which \cite{k,mcs, V1}. Targets with a linear structure have continued to spark research interest, see e.g.~\cite{jl, jls, ln}. Target spaces which have non-positive sectional curvature were considered in \cite{V2,ls,lps,ln}, see also the references therein.
In his seminal paper \cite{Gromov-cc}, Gromov considered Carnot-Carath\'eodory spaces as target spaces. He proved in particular that every partially defined Lipschitz map from a compact Riemannian manifold $M$ of dimension $n$ to a contact manifold of dimension $\geq 2n+1$, endowed with a Carnot-Carath\'eodory metric, has a Lipschitz extension to all of $M$. Gromov's proof uses microflexibility. Later, Allcock gave a more elementary proof of the special case that a map from the circle to a higher Heisenberg group can be extended to the disc \cite{Allcock}.  A detailed or more elementary proof of Gromov's full result which is not based on microflexibility does not seem to exist in the literature.

The aim of this note is to provide such an elementary and detailed proof of Gromov's result in the case when the target space is a jet space Carnot group $J^k(\R^n)$; see below for definitions. These spaces, which give a model for a certain class of $(k+1)$-step Carnot groups including the $n$-th Heisenberg group $J^1(\R^n)$, have recently been considered by the second author \cite{Young-nilpotent-isop} in connection with isoperimetric inequalities for higher-dimensional cycles.  The methods therein circumvent the use of microflexibility with a construction involving jet maps.
We will combine some of the ideas from \cite{Young-nilpotent-isop} with techniques recently developed
by Lang-Schlichenmaier \cite{Lang-Schlichenmaier} in order to prove Lipschitz extension theorems for partially defined Lipschitz maps from a rather general domain space $X$ to $J^k(\R^n)$, which includes the case when $X$ is a compact Riemannian manifold of dimension $\leq n$.

Before we state our results, recall that a Carnot group is a connected, simply connected Lie group with stratified Lie algebra. The left-invariant sub-Riemannian Carnot-Carath\'eodory distance $d_c$ is defined with respect to the first, so-called horizontal, layer of the stratification. The jet spaces $J^k(\R^n)$ give a model for a certain class of $(k+1)$-step Carnot groups, including the $n$-th Heisenberg group $J^1(\R^n)$, the Engel group $J^2(\R)$ and more generally the model filiform groups $J^k(\R)$. Definitions will be given in \secref{Section:prelims}. A special case of our main result can be stated as follows:
\bt\label{thm:riemannian}
 Let $n, k\geq 1$ and let $M$ be a compact Riemannian manifold, possibly with boundary, of dimension at most $n$. Then there exists a constant $C$ depending only on $n$ and $M$ such that for every nonempty subset $Z\subset M$, every $\lambda$-Lipschitz map $f: Z\to\ (J^k(\R^n), d_c)$ has a $C\lambda$-Lipschitz extension $\bar{f}: M\to (J^k(\R^n), d_c)$.
\et

Our theorem also holds for more general metric spaces.  In the following, $\ANdim(X)$ denotes the Assouad-Nagata dimension of the metric space $X$.  Its definition will be given in \secref{Section:prelims}.
\bt\label{thm:mainCombined}
 Let $X$ be a metric space and $n, k\geq 1$. Let $Z\subset X$ be a nonempty closed subset.  Suppose either $\ANdim(Z)\le n-1$ or $\ANdim(X\backslash Z)\le n$. Then there exists a constant $C$ depending only on $n$ and on the implicit constant in $\ANdim$ such that every $\lambda$-Lipschitz map $f: Z\to\ (J^k(\R^n), d_c)$ has a $C\lambda$-Lipschitz extension $\bar{f}: X\to (J^k(\R^n), d_c)$.
\et

We mention here that every subset of $\R^n$ with nonempty interior  has Assouad-Nagata dimension $n$, and so does every compact metric space locally bilipschitz homeomorphic to such subsets. In particular, every compact Riemannian $n$-manifold (with boundary) has Assouad-Nagata dimension $n$, so Theorem \ref{thm:riemannian} is a corollary of Theorem \ref{thm:mainCombined}. Furthermore, homogeneous Hadamard $n$-manifolds with pinched curvature, products of $n$ trees, and Euclidean buildings of rank $n$ all have Assouad-Nagata dimension $n$, see \cite{Lang-Schlichenmaier}. 

In general, Lipschitz maps from the sphere $S^n$ to $(J^k(\R^n), d_c)$ need not admit Lipschitz extensions to the unit ball $\bar{B}^{n+1}$ as was shown in \cite{Rigot-Wenger}, see also \cite{Balogh-Faessler} for the case $k=1$.

A result analogous to Theorem \ref{thm:mainCombined} with $(J^k(\R^n), d_c)$ replaced by a Lipschitz $(n-1)$-connected metric space $Y$ is due to Lang-Schlichenmaier \cite{Lang-Schlichenmaier}. Recall that a space $Y$ is Lipschitz $(n-1)$-connected if there exists $c$ such that every $\lambda$-Lipschitz map $f: S^m\to Y$, with $m\leq n-1$, possesses a $c\lambda$-Lipschitz extension $\bar{f}: \bar{B}^{m+1}\to Y$. In order to obtain our results above, it would therefore be enough to prove that $(J^k(\R^n), d_c)$ is Lipschitz $(n-1)$-connected. It turns out, however, that it is about as difficult to construct Lipschitz extensions from $Z$, where $Z$ is as in our theorems, as it is from $S^m$.

The proof of \thmref{thm:mainCombined} proceeds along the lines of the proof of the analogous result by Lang-Schlichenmaier mentioned above. Roughly, the idea is to first construct a suitable (abstract) simplicial complex $\Sigma$ of dimension at most $n$ and a map $g: X\backslash Z\to \Sigma$. In some sense, the pair $(\Sigma, g)$ provides a decomposition of $X\backslash Z$ into simplices of diameter proportional to the distance to $Z$. In \cite{Lang-Schlichenmaier} the authors then construct a map $h: \Sigma\to Y$ by first choosing for each vertex $e_i$ of $\Sigma$ a suitable point $z_i$ in $Z$, such as a nearest neighbor in $Z$, then setting $h(e_i):= f(z_i)$ and finally extending $h$ to higher skeleta of $\Sigma$ and to all of $\Sigma$ by means of the Lipschitz connectedness of $Y$. The Lipschitz extension $\bar{f}$ of $f$ is obtained by setting $\bar{f} = h\circ g$ on $X\backslash Z$.

In our situation, in which the Lipschitz connectedness is not available, we combine and adapt ideas and arguments from \cite{Lang-Schlichenmaier} and \cite{Young-nilpotent-isop} to construct a suitable map $h: \Sigma\to J^k(\R^n)$. 
We first remark that one has an abundance of Lipschitz maps from simplicial complexes to
$J^k(\R^n)$.  For instance, any map from the $0$-skeleton
$\Delta^{(0)}$ of the standard $n$-simplex $\Delta$ to $J^k(\R^n)$ can
easily be extended to a Lipschitz map from $\Delta$ if the images of
the vertices lie ``in general position" \cite{Young-nilpotent-isop},
but the Lipschitz constant depends on the vertices.  
We will use a technique from \cite{Young-nilpotent-isop} to bound the Lipschitz constants of these extensions.  Using the group structure and the scaling automorphisms of $J^k(\R^n)$, we can rescale, perturb, and translate a map from the $0$-skeleton of an $n$-simplex so that its image lies in a fixed finite set.  Such
maps can be extended using a finite family of extensions, giving us an automatic bound on the Lipschitz
constant.  When we apply this construction to $\Sigma$, the necessary
perturbations are small for simplices close to $Z$, so the
construction produces a Lipschitz extension.

Gromov used
microflexibility to construct extensions from $\Delta^{(0)}$ to $\Delta$
continuously with respect to the vertices.  This gives a uniform bound
on the Lipschitz constant as long as the vertices stay in a compact
set, and any
set of vertices may be rescaled to lie in a compact set.  This implies
a Lipschitz extension theorem for the pair $(\Delta^{(0)},\Delta)$.
Extending this result to a complex $\Sigma$ takes some additional
care to ensure that extensions on adjacent simplices agree on their
intersection.

\bigskip

{\noindent \bf Acknowledgements:} The second author would like to thank the Courant Institute of Mathematical Sciences in New York for its hospitality during part of the writing of this paper.

\section{Preliminaries}\label{Section:prelims}
\subsection{Nagata dimension}
We recall here the definition and some properties of the Assouad-Nagata dimension for a metric space $Y$. For a detailed account we refer to \cite{Lang-Schlichenmaier}.  A family $(B_i)_{i\in I}$ of subsets of $Y$ is called $D$-bounded if $\diam B_i\leq D$ for all $i\in I$. For $s\geq 0$, the $s$-multiplicity of the family is the infimum of all $k\geq 0$ such that every subset of $X$ with diameter $\leq s$ intersects at most $k$ members of the family. 

\bd
 The (Assouad-)Nagata dimension $\ANdim(X)$ of a metric space $Y$ is the infimum of all integers $n$ with the following properties: there exists a constant $c>0$ such that for all $s>0$, $Y$ has a $cs$-bounded covering with $s$-multiplicity at most $n+1$.
\ed

It can be shown that the Nagata dimension of $Y$ is at least its topological dimension. If $Y=Y_1\cup Y_2$ then the Nagata dimension of $Y$ is the maximum of the Nagata dimensions of $Y_1$ and $Y_2$. Every subset of $\R^n$ with nonempty interior has Nagata dimension $n$. The Nagata dimension is invariant under biLipschitz homeomorphisms. It thus follows that every compact metric space locally biLipschitz homeomorphic to an open subset of $\R^n$ has Nagata dimension $n$. In particular, every compact Riemannian $n$-manifold has Nagata dimension $n$. For proofs of all these statements and many more properties see \cite{Lang-Schlichenmaier}. 

In the proof of our main theorems we will need the following technical result of Lang-Schlichenmaier \cite{Lang-Schlichenmaier}. 

\bt\label{thm:good-covering}
 Let $X$ be a metric space and $Z\subset X$ a nonempty closed set. Suppose that either 
 $\ANdim(Z)\leq n-1$ for some $n\geq 1$ or $\ANdim(X\backslash Z)\leq n$ for some $n\geq 0$, with constant $c$. Then there exist $\alpha, \beta > 0$, only depending on $n$ and $c$, and  a covering $(B_i)_{i\in I}$ of $X\backslash Z$ by subsets of $X\backslash Z$ such that
\begin{enumerate}
 \item $\diam B_i \leq \alpha d(B_i, Z)$ for every $i\in I$,
 \item every set $D\subset X\backslash Z$ with $\diam D\leq \beta d(D,Z)$ meets at most $n + 1$ members of $(B_i)_{i\in I}$.
\end{enumerate}
\et

In the case $\ANdim(X\backslash Z)\leq n$ the existence of such a covering $(B_i)_{i\in I}$ is established in the proof of Theorem 1.5 in \cite{Lang-Schlichenmaier}. In the case $\ANdim(Z)\leq n-1$ the existence is established in the proof of Theorem 1.6 in \cite{Lang-Schlichenmaier}. \thmref{thm:good-covering} will allow us to construct suitable simplicial complexes which will be used in the construction of a Lipschitz extension.

\subsection{Jet spaces as Carnot groups}
The jet spaces $J^k(\R^n)$ give a model for a certain class of Carnot groups, which includes a model for the Heisenberg group $J^1(\R^n)$, for the Engel group $J^2(\R)$ and more generally for the model filiform groups $J^k(\R)$. We refer to \cite{Warhurst} for a detailed account on these groups.

Let $n, k\geq 1$. The $k$-th order Taylor polynomial of a $C^k$-smooth function $f:\R^n \rightarrow \R$ at $x_0$ is given by
\begin{equation*}
 T^k_{x_0}(f) (x) = \sum_{j=0}^k \sum_{I\in I(j)} \p_I f(x_0) \, \dfrac{(x-x_0)^I}{I!},
\end{equation*}
where $I(j)$ denotes the set of all $j$-indices, that is, of $n$-tuples $I = (i_1,\dots,i_n)$ satisfying $|I| = i_1 + \cdots + i_n = j$, and where $I!=i_1!\dots i_n!$,  furthermore $x^I=(x_1)^{i_1}\cdots (x_n)^{i_n}$, and 
\begin{equation*}
 \p_I f(x_0) = \dfrac{\p^j f}{\p x_1^{i_1}\dots \p x_n^{i_n}} (x_0).
\end{equation*}
Two functions $f_1$, $f_2 \in C^k(\R^n,\R)$ are said to be equivalent at $x_0$, $f_1 \sim_{x_0} f_2$, if and only if $T^k_{x_0}(f_1) =  T^k_{x_0}(f_2)$. The equivalence class of $f$ at $x_0$ is denoted by $j_{x_0}^k(f)$. The $k$-jet space over $\R^n$ is 
\begin{equation*}
 J^k(\R^n) = \bigcup_{x_0\in\R^n} \, C^k(\R^n,\R)/\sim_{x_0},
\end{equation*}
endowed with the topology such that
\begin{equation*}
 J^k(\R^n) \equiv \R^n \times \R^{d_k^n} \times \dots \times \R^{d_0^n},\quad\text{where $d_j^n = \binom{n+j-1}{j}$,}
\end{equation*}
identified via the global coordinates $(x,u^k,\dots,u^0)$, where $x(j_{x_0}^k(f)) = x_0$ and $u^j = (u_I^j)_{I\in I(j)}$ with $u_I^j(j_{x_0}^k(f)) = \p_I f(x_0)$ for $j=0,\dots,k$. For an element $a\in J^k(\R^n)$, given in global coordinates, we will often write $\pi(a)$ to denote the $x$-coordinate of $a$ and $a_I$, for $I\in I(j)$, to denote the $u_I^j$-coordinate of $a$.   

The so-called horizontal subbundle of the tangent bundle $TJ^k(\R^n)$ of $J^k(\R^n)$ is defined by
\begin{equation*}
\mathcal{H}_0 = \spa \left\{X_i:  i=1,\dots,n\right\} \oplus \spa \left\{\p_{u_I^k}: I \in I(k)\right\},
\end{equation*}
where
\begin{equation*}
X_i = \p_{x_i} + \sum_{j=0}^{k-1} \sum_{I\in I(j)} u_{I+e_i}^{j+1} \p_{u_I^j}
\end{equation*}
with $e_i=(0,\dots, 1, \dots, 0) \in I(1)$.
For $j=1,\dots,k$, we set
\begin{equation*}
\mathcal{H}_j = \spa \left\{\p_{u_I^{k-j}}: I \in I(k-j)\right\}.
\end{equation*}
The only non trivial commutators are 
\begin{equation*}
 \left[\p_{u_{I+e_i}^{j+1}}, X_i\right] = \p_{u_I^j}, \quad I\in I(j), \quad j=0,\dots,k-1.
\end{equation*}
It follows that $\mathcal{H}_j = [\mathcal{H}_0,\mathcal{H}_{j-1}]$ for $j=1,\dots,k$ and $\mathcal{H} = \mathcal{H}_0 \oplus \dots \oplus \mathcal{H}_k$ is a $(k+1)$-step stratified Lie algebra which spans $TJ^k(\R^n)$ pointwise. Corresponding to this Lie algebra there is a Carnot group $G^k_n$, i.e., a connected, simply connected and stratified Lie group, unique up to isomorphism. As is shown in \cite{Warhurst}, the product $\odot$ which makes $(J^k(\R^n), \odot)$ into a Carnot group isomorphic to $G^k_n$, is given by
\begin{equation}\label{eqn:product}
 j^k_x(f) \odot j^k_y(g):= j^k_{x+y}(h)\quad\text{ where }\quad  h(z):= T_x^kf (z) + T^k_yg(z-x).
\end{equation}
For points $a$ and $b$, given in global coordinates, this becomes
\begin{equation*}
 \pi(a\odot b) = \pi(a) + \pi(b)
\end{equation*}
and
\begin{equation*}
 (a\odot b)_I = b_I + \sum_{J\geq I} a_J \frac{\pi(b)^{J-I}}{(J-I)!},
\end{equation*}
where, for $J=(j_1, \dots, j_n)$ and $I=(i_1, \dots, i_n)$, we write $J\geq I$ if and only if $j_m\geq i_m$ for all $m$. We will often write $a\cdot b$ or even $ab$ instead of $a\odot b$.

The dilation homomorphisms $\delta_L$ on $J^k(\R^n)$ are defined by
\begin{equation*}
 \delta_L(j^k_x(f)) := j^k_{Lx}(f_L)\quad\text{ where }\quad f_L(z):= L^{k+1}f(z/L).
\end{equation*}
In global coordinates this becomes
\begin{equation*}
 \delta_L(x,u^k,u^{k-1},\dots,u^0) = (Lx,Lu^k,L^2 u^{k-1}, \dots,L^{k+1} u^0).
\end{equation*}

In the next section, we will need the following constructions. Given $a\in J^k(\R^n)$ and a smooth function $f: V\to\R$, where $V$ is an open subset of $\R^n$, we define the smooth function $f^a: V + \pi(a) \to \R$ by
\begin{equation*}
 f^a(z):= f(z-\pi(a)) + \sum_Ia_I\frac{(z-\pi(a))^I}{I!}.
\end{equation*}
It is not difficult to check that $(f^a)^b = f^{b\odot a}$ for all $a,b\in J^k(\R^n)$ and that $$j^k_{x+\pi(a)}(f^a) = a\odot j^k_x(f)$$ for every $x\in V$.

\subsubsection*{The Carnot-Carath\'eodory distance on $J^k(\R^n)$}
Let $g_0$ be the left invariant Riemannian metric such that $(X_1,\dots,X_n,\p_{u_I^j})_{j=0,\dots,k, I\in I(j)}$ is an orthonormal basis.

The Carnot-Carath\'eodory distance $d_c$ on $J^k(\R^n)$ is the sub-Riemannian distance defined by 
\begin{equation*}
 d_c(x,y) = \inf \{ length_{g_0} (\gamma); \; \gamma \text{ horizontal } C^1 \text{ curve joining } x \text{ to } y\},
\end{equation*}
where a $C^1$ curve is said to be horizontal if, at every point, its tangent vector belongs to the horizontal subbundle of the tangent bundle. The Carnot-Carath\'eodory distance is left invariant and moreover $1$-homogeneous with respect to the dilations, i.e., $d_c(\delta_L(x),\delta_L(y)) = L\,d_c(x,y)$ for all $x$, $y \in J^k(\R^n)$ and all $L\geq 0$. 

\subsubsection*{Jet maps}
A differentiable map from an open subset $U\subset\R^n$ with values in $J^k(\R^n)$ is said to be horizontal if the image of its differential lies is the horizontal subbundle of the tangent bundle. If $f\in C^{k+1}(U,\R)$, then the map $j^k(f):x \mapsto j_{x}^k(f)$ is a horizontal $C^1$ map with 
\begin{equation*}
 \p_{x_i} (j^k(f)) (x) =  X_i(j_x^k(f)) + \sum_{I \in I(k)} \p_{I+e_i} f(x) \,\p_{u_I^k}
\end{equation*}
for $i=1,\dots,n$. Fix $x,y\in\R^n$ and define $\gamma(t) := (1-t)x+ty$. It follows that $t\in [0,1] \mapsto j^k(f)(\gamma(t))$ is a horizontal curve between $x$ and $y$ and we get from the very definition of $d_c$ that 
\begin{equation*}
 d_c(j^k(f)(x),j^k(f)(y)) \leq \sup_{t\in [0,1]} \left(1+\sum_{I\in I(k)} \sum_{i=1}^n (\p_{I+e_i} f(\gamma(t))^2\right)^{1/2} \, \|y-x\|.
\end{equation*}
In particular, the map $j^k(f):\R^n \rightarrow (J^k(\R^n),d_c)$ is locally Lipschitz.

\section{Lipschitz maps from simplices}

Given a nonempty set $A$ we define the simplex $\Sigma(A)$ by
\begin{equation*}
 \Sigma(A):=  \left\{(v_a)\in l^2(A): v_a\geq 0, \sum_{a\in A} v_a=1\right\}
\end{equation*}
and denote by $\Sigma^{(n)}(A)$ its $n$-skeleton for $n\in\N\cup\{0\}$.
If $A$ is a subset of a metric space $(Y, d)$ and $\varepsilon>0$ we furthermore set
\begin{equation*}
 \Sigma(A, \varepsilon, n):= \left\{[e_{a_1}, \dots, e_{a_m}]\subset\Sigma^{(n)}(A): \text{ $m\leq n+1$, $d(a_i, a_j)\leq \varepsilon$ for $1\leq i,j\leq m$}\right\},
\end{equation*}
where $e_a\in\Sigma^{(0)}(A)$ denotes the vertex corresponding to the index $a$ and where $[e_{a_1}, \dots, e_{a_m}]$ is the simplex spanned by the vectors $e_{a_1}, \dots, e_{a_m}$. Each simplex comes with the restricted metric from $l^2(A)$.

The following proposition is a variation of Lemma 12 in \cite{Young-nilpotent-isop}.

\bp\label{prop:good-subset-Heis}
 Let $k, n\geq 1$ and $r\geq 2$ be integers, and let $\varepsilon\geq \varepsilon'>0$. Then there exist discrete subsets $A, A'\subset J^k(\R^n)$, each $\varepsilon'$-dense with respect to the $d_c$-distance, and maps 
 $$F: \Sigma(\Gamma, \varepsilon, n) \to J^k(\R^n)\quad\text{ and }\quad F': \Sigma(\Gamma', \varepsilon, n) \to J^k(\R^n),$$
 where $\Gamma= A\cup \delta_r(A')$ and $\Gamma'=A'\cup \delta_r(A)$,
 such that the following properties hold: Firstly, $F(e_c) = c$ and $F'(e_{c'}) = c'$ for all $c\in \Gamma$ and $c'\in \Gamma'$; secondly, for $m\leq n+1$ and $a_1, \dots, a_m\in A$ with $d_c(a_i, a_j)\leq \varepsilon/r$ 
   \begin{equation}\label{eq:F-prop}
    F'{|_{[e_{\delta_r(a_1)}, \dots, e_{\delta_r(a_m)}]}} = \delta_r\circ  F{|_{[e_{a_1}, \dots, e_{a_m}]}}
   \end{equation}
   and, similarly, for $a'_1, \dots, a'_m\in A'$ with $d_c(a'_i, a'_j)\leq \varepsilon/r$
   \begin{equation}\label{eq:F'-prop}
     F{|_{[e_{\delta_r(a'_1)}, \dots, e_{\delta_r(a'_m)}]}} = \delta_r\circ F'{|_{[e_{a'_1}, \dots, e_{a'_m}]}};
   \end{equation}
 thirdly, the restrictions of $F$ and $F'$ to any simplex are $\varrho$-Lipschitz with respect to the $d_c$-distance on $J^k(\R^n)$, for some $\varrho$ independent of the simplex.
\ep

The simplex appearing on the left hand side of \eqref{eq:F-prop} should be thought as identified with the simplex appearing on the right hand side via the natural isometry which sends $e_{\delta_r(a_i)}$ to $e_{a_i}$. The same remark applies to \eqref{eq:F'-prop}.

\begin{proof}
Our first step is to construct $A$ and $A'$.  We want $A$ and $A'$ to satisfy two properties.  First,
if $[e_{v_0},\dots, e_{v_m}]$ is a simplex in $\Sigma(\Gamma, \varepsilon, n)$ or $\Sigma(\Gamma', \varepsilon, n)$, then  $\pi(v_0), \dots, \pi(v_m)$ span a non-degenerate $m$-simplex in $\R^n$.  Second, we want the set of isometry classes of the simplices $[\pi(v_0), \dots, \pi(v_m)]$ to be finite; we will accomplish this goal by choosing $A$ and $A'$ to be $\Lambda$-invariant, where $\Lambda$ is a lattice in $J^k(\R^n)$.

 Let $\Lambda\subset J^k(\R^n)$ be a cocompact lattice such that $\delta_L(\Lambda)\subset\Lambda$ for all $L\in\N$, see e.g.~Proposition 2 in \cite{Young-nilpotent-isop}.  Let $R>0$ be sufficiently large so that the closed $d_c$-ball $B:=\bar{B}(0,R)$ satisfies $\Lambda B=J^k(\R^n)$. Let $\eta>0$ be small enough so that $d_c(x,\gamma\cdot x)>\eta$ for all $x\in B$ and $\gamma\in\Gamma\backslash\{0\}$. After rescaling the metric $d_c$ we may assume that $\varepsilon<\eta$. 
 In the following, a subset $C\subset J^k(\R^n)$ will be called admissible if for all points $v_0, \dots, v_m\in C\cap\bar{B}(0, rR+ r\eta)$ satisfying $m\leq n$ and
 \begin{equation}\label{eqn:distance-vi-vj}
 0< d_c(v_i, v_j) < \varepsilon \quad\text{for all $i\not=j$,}
 \end{equation}
 $\pi(v_0), \dots, \pi(v_m)$ span a non-degenerate $m$-simplex in $\R^n$. Note that if $C$ is admissible then so is $\delta_r(C)$.  

Let $\{y_1,\dots, y_t\}\subset B$ be a finite and $\varepsilon'/2$-dense subset of $B$. Choose $x_1, x'_1,\dots, x_t, x'_t\in B$ successively in such a way that, for all $i=1, \dots, t$, 
\begin{itemize}
\item $d_c(x_i, y_i)<\varepsilon'/2$ and $d_c(x'_i, y_i)<\varepsilon'/2$,
\item $x_i\not\in \Lambda\{x_1,\dots, x_{i-1}\}$ and $x'_i\not\in \Lambda\{x'_1,\dots, x'_{i-1}\}$ when $i\geq 2$,
\item $\Lambda\{x_1, \dots, x_i\}$ and $\delta_r(\Lambda\{x'_1,\dots, x'_i\})$ are disjoint and their union is admissible,
\item $\Lambda\{x'_1, \dots, x'_i\}$ and $\delta_r(\Lambda\{x_1,\dots, x_i\})$ are disjoint and their union is admissible.
\end{itemize}
Set $A:= \Lambda\{x_1,\dots, x_t\}$ and $A':= \Lambda\{x'_1,\dots, x'_t\}$. The sets $A$ and $A'$ are $\varepsilon'$-dense in $J^k(\R^n)$. Set $\Gamma= A\cup \delta_r(A')$ and $\Gamma'=A'\cup \delta_r(A)$; let $\mathscr{G}$ and $\mathscr{G}'$ denote the collections of simplices in $\Sigma(\Gamma, \varepsilon, n)$ and $\Sigma(\Gamma', \varepsilon, n)$, respectively. For each $(m-1)$-simplex $S= [e_{b_1}, \dots, e_{b_m}]\in\mathscr{G}\cup\mathscr{G}'$ denote by $\bar{\pi}(S)$ the closed convex hull in $\R^n$ of the points $\pi(b_1), \dots, \pi(b_m)$. Note that $\{\bar{\pi}(S): S\in\mathscr{G}\cup\mathscr{G}'\}$ contains only finitely many isometry classes of simplices. Let $\varepsilon_1>\varepsilon_2>\dots>\varepsilon_{n+1}>0$ be suitably small, to be determined later, and define for each $S=[e_{b_1}, \dots, e_{b_m}]\in\mathscr{G}\cup\mathscr{G}'$ an open subset $U_S\subset \R^n$ as follows. If all $b_i$ are in $\delta_r(A)$ or all $b_i$ are in $\delta_r(A')$ then let $U_S\subset\R^n$ be the open $r\varepsilon_m$-neighborhood of $\bar{\pi}(S)$; otherwise, let $U_S\subset\R^n$ be the open $\varepsilon_m$-neighborhood of $\bar{\pi}(S)$. Clearly, if the $\varepsilon_m$ were chosen suitably (independently of $S$) then for distinct faces $S, T$ of the same simplex we have
\begin{equation*}
 U_S\cap U_T\subset \frac{1}{2}U_{S\cap T},
\end{equation*}
where we set $U_{S\cap T}:= \emptyset$ if $S\cap T = \emptyset$, and where $\frac{1}{2}U_{S\cap T}$ denotes the neighborhood of half the radius of that of $U_{S\cap T}$. Note that if $S=[e_{a_1}, \dots, e_{a_m}]$ is a simplex in $\mathscr{G}\cup\mathscr{G}'$ and $\gamma$ is such that $\gamma S:=[e_{\gamma a_1}, \dots, e_{\gamma a_m}]$ is also a simplex in $\mathscr{G}\cup\mathscr{G}'$ then $U_{\gamma S} = \pi(\gamma) + U_S$.


Next, we associate with each $S\in\mathscr{G\cup\mathscr{G}'}$ a
smooth function $f_S:U_S\to\R$ such that if $S\in
\mathscr{G\cup\mathscr{G}'}$ and $S'$ is a face of $S$, then $f_S$ and
$f_{S'}$ coincide on $U_{S}\cap \frac{1}{2} U_{S'}$.  In particular,
if $S$ and $T$ are simplices which share a face, then they coincide on
$U_{S}\cap U_T \cap \frac{1}{2} U_{S\cap T}$.
We first construct the $f_S$ for simplices in $\Sigma(A, \varepsilon,
n)$.  For this, let $\mathscr{A}$ be the collection of all simplices
in $\Sigma(A, \varepsilon, n)$ and
$\mathscr{A}^{(d)}\subset\mathscr{A}$ the subcollection of 
$d$-simplices.  For $S= [e_{x_i}]$ with $1\leq i\leq t$ choose $f_S$
with $j^k_{\pi(x_i)}(f_S) =x_i$, and for each $\gamma\in\Lambda$, let
$f_{\gamma S}:= (f_S)^\gamma$.  This yields a function $f_S$ for each
$S\in\mathscr{A}^{(0)}$.  Now, we will use induction on $d$ to define
functions $f_S$ for each $S\in \mathscr{A}^{(d)}$.
Note that $\Lambda$ acts on $\mathscr{A}^{(d)}$ and there are finitely
many orbits.  We will define the $f_S$ to be $\Lambda$-equivariant;
that is, if $\gamma\in \Lambda$ and $S\in
  \mathscr{A}$,  then $(f_{S})^\gamma=f_{\gamma S}$.
Let $S_1,\dots, S_q$ be representatives of the $\Lambda$-orbits.
If $1\le i \le q$, we may assume by induction that $f_{T}$ is defined
for each face $T$ in the boundary of $S_i$.  Furthermore, if $T$ and
$T'$ are both codimension-1 faces of $S_i$, then $U_T\cap U_{T'}\subset \frac{1}{2} U_{T\cap T'},$
so by induction,
$$f_T|_{U_T\cap U_{T'}} = f_{T\cap T'}|_{U_T\cap U_{T'}} = f_{T'}|_{U_T\cap U_{T'}}.$$
This allows us to define a smooth function $f'$ on $\bigcup_T
\frac{1}{2} U_T$, where the union is taken over codimension-1 faces of $S_i$.  We can extend $f'$ to a smooth function $f_{S_i}$ defined on $U_{S_i}$.
If $S\in \mathscr{A}^{(d)}$, then $S=\gamma S_i$ for some $\gamma\in \Lambda$ and some $1\le i\le q$.  Let $f_S=(f_{S_i})^\gamma$.
Now let $\mathscr{A}_r$ denote the collection of simplices in $\Sigma(\delta_r(A), \varepsilon, n)$. With each $S=[e_{\delta_r(a_1)}, \dots, e_{\delta(a_m)}]\in\mathscr{A}_r$ we associate the function $f_S: U_S\to\R$ defined by $f_S(z):= r^{k+1}f_{[e_{a_1},\dots, e_{a_m}]}(z/r)$. 

Let $\mathscr{A}'$ and $\mathscr{A}'_r$ be the collections of
simplices in $\Sigma(A', \varepsilon, n)$ and $\Sigma(\delta_r(A'),
\varepsilon, n)$, respectively. We may proceed exactly as above to
associate with every simplex $S\in\mathscr{A}'\cup \mathscr{A}'_r$ a
smooth function $f_S: U_S\to\R$. Finally, consider $\mathscr{G}$.  We
have associated functions $f_S$ to the subset $\mathscr{A}\cup
\mathscr{A}'_r\subset \mathscr{G}$, and we can proceed again as above
to obtain a function $f_S$ for every $S\in\mathscr{G}$. The main
difference is that $\delta_r(\Lambda)\subset \Lambda$, rather than
$\Lambda$, acts on $\mathscr{G}$, so we choose the $S_i$ to be
representatives of the $\delta_r(\Lambda)$-orbits of $\mathscr{G}\setminus(\mathscr{A}\cup
\mathscr{A}'_r)$.  In the same way we obtain functions $f_S$ for $S$ in the collection $\mathscr{G}'$ of simplices in $\Sigma(\Gamma', \varepsilon, n)$.

Finally, for each element $S= [e_{b_1},\dots, e_{b_m}]$ in $\mathscr{G}\cup\mathscr{G}'$ define a map $\bar{\pi}_S:S\to\R^n$ by $\bar{\pi}_S((v_i)_i):= v_1\pi(b_1) + \dots v_m\pi(b_m)$. Note that each $\bar{\pi}_S$ is bilipschitz continuous with a bilipschitz constant independent of $S$ because  $\{\bar{\pi}(S): S\in\mathscr{G}\cup\mathscr{G}'\}$ contains only finitely many isometry classes of simplices. Let $F: \Sigma(\Gamma, \varepsilon, n) \to J^k(\R^n)$ and $F': \Sigma(\Gamma', \varepsilon, n) \to J^k(\R^n)$ be the maps such that $F|_S= j^k(f_S)\circ \bar{\pi}_S$ for each $S\in\mathscr{G}$, and analogously  $F'|_S= j^k(f_S)\circ \bar{\pi}_S$ for each $S\in\mathscr{G}'$. These maps are well-defined because the functions $f_S$ agree on neighborhoods of common faces of simplices. Furthermore, we have
\begin{equation*}
 F(e_c) = j^k_{\bar{\pi}_{[e_c]}(e_c)}(f_{[e_c]}) = c
\end{equation*}
for all $c\in\Gamma$, and analogously for $F'$. We can also compute for $S=[e_{a_1}, \dots, e_{a_m}]\in\mathscr{A}$ and $(v_i)\in S$ 
\begin{equation*}
 \delta_r\circ  F{|_S}((v_i)) = \delta_r(j^k_{\bar{\pi}_S((v_i))}(f_S)) = j^k_{r\bar{\pi}_S((v_i))}(r^{k+1}f_S(\cdot /r)) = F'{|_{[e_{\delta_r(a_1)}, \dots, e_{\delta_r(a_m)}]}}((v_i)),
\end{equation*}
from which the second property in the statement of the proposition follows.
Finally, it is straightforward to check that  $F|_{\gamma S} = \gamma\odot F|_S$ and $F'|_{\gamma S} = \gamma\odot F'|_S$. The Lipschitz property is a consequence of this and of the fact that $j^k(f_S)$ is Lipschitz. This concludes the proof.
\end{proof}

\section{Proof of the main results}

Theorem \ref{thm:mainCombined} is a direct consequence of the following theorem together with \thmref{thm:good-covering}.

\bt\label{thm:technical-extension}
 Let $X$ be a metric space, and let $\alpha, \beta > 0$ and $n,k\geq 1$. Suppose that
$Z\subset X$ is a nonempty closed set and $(B_i)_{i\in I}$ is a covering of $X\backslash Z$ by subsets of $X\backslash Z$ such that
\begin{enumerate}
 \item $\diam B_i \leq \alpha d(B_i, Z)$ for every $i\in I$,
 \item every set $D\subset X\backslash Z$ with $\diam D\leq \beta d(D,Z)$ meets at most $n + 1$ members of $(B_i)_{i\in I}$.
\end{enumerate}
Then there exists a constant $C$ depending only on $\alpha, \beta, n$ such that every $\lambda$-Lipschitz map $f: Z\to (J^k(\R^n), d_c)$ has a $C\lambda$-Lipschitz extension $\bar{f}: X\to (J^k(\R^n), d_c)$.
\et

The proof of \thmref{thm:technical-extension} follows in parts the proof of Theorem 5.2 in \cite{Lang-Schlichenmaier}.

\begin{proof}
 In the following, we will write $J^k(\R^n)$ to mean $J^k(\R^n)$ endowed with the metric $d_c$. Let $f:Z\to J^k(\R^n)$ be a $\lambda$-Lipschitz map. After rescaling the metric on $X$ by the factor $1/\lambda$ we may assume 
 that $\lambda=1$, i.e.~ $f$ is $1$-Lipschitz. Let $(B_i)_{i\in I}$ be as in the 
 hypothesis. Of course, we may assume that the sets $B_i$ are all nonempty. Set $\tau:=\beta/(2(\beta+1))$ and note that $\tau< 1/2$. For each $i\in I$, define a $1$-Lipschitz function $\sigma_i: X\backslash Z\to\R$ by
\begin{equation*}
 \sigma_i(x):= \max\{0, \tau d(B_i, Z) - d(x, B_i)\}.
\end{equation*}
It can be shown exactly as in \cite{Lang-Schlichenmaier} that for every $x\in X\backslash Z$ there are at most $n+1$ indices $i\in I$ with $\sigma_i(x)>0$. Indeed, for every such $i$, pick $x_i\in B_i$ with $d(x, x_i)<\tau d(B_i, Z)\leq \tau d(x_i, Z)$. The set $D$ of all these $x_i$ satisfies $\diam D\leq 2\tau \sup_i d(x_i, Z)\leq 2\tau(\diam D + d(D, Z))$, thus $\diam D\leq \beta d(D, Z)$. By condition (ii), $D$ therefore meets at most $n+1$ members of $(B_i)_{i\in I}$, thus $\sigma_i(x)>0$ for at most $n+1$ indices $i\in I$, as claimed. Now, let $\bar{\sigma}:= \sum_{i\in I}\sigma_i$ and note that $\bar{\sigma}>0$ on $X\backslash Z$. Define $g: X\backslash Z\to l^2(I)$ by 
\begin{equation*}
 g(x):= \left(\frac{\sigma_i(x)}{\bar{\sigma}(x)}\right)_{i\in I}.
\end{equation*}
Set $\Sigma:= \Sigma(I)$
and observe that
\begin{equation*}
 g(X\backslash Z) \subset \Sigma':= \left\{[e_{i_1}, \dots, e_{i_m}] \subset \Sigma^{(n)}: \text{ $\exists x \in X\backslash Z$ with $\sigma_{i_j}(x)>0$ for $j=1, \dots, m$}\right\}.
\end{equation*}
Here, $e_i\in \Sigma^{(0)}$ denotes the vertex corresponding to the index $i\in I$. Let $x\in X\backslash Z$ and $i\in I$ with $\sigma_i(x)>0$. We claim that 
\begin{equation}\label{eqn:x-B-Z}
 (1+\alpha + \tau)^{-1} d(x, Z)< d(B_i, Z) < (1-\tau)^{-1} d(x, Z).
\end{equation}
Indeed, (i) and the fact that $d(x, B_i)< \tau d(B_i, Z)$ yield
\begin{equation*}
 d(x, Z) \leq d(x, B_i) + \diam B_i + d(B_i, Z) < (1 + \alpha + \tau)d(B_i, Z)
\end{equation*}
and
\begin{equation*}
 d(B_i, Z) \leq d(x, B_i) + d(x, Z) < \tau d(B_i, Z) + d(x, Z),
 \end{equation*}
 thus both inequalities in \eqref{eqn:x-B-Z}.
Next, for each $i\in I$ choose $z_i\in Z$ satisfying
\begin{equation*}
 d(z_i, B_i) \leq (2-\tau)d(B_i, Z)
\end{equation*}
and note that for $x\in X\backslash Z$ and $i\in I$ with $\sigma_i(x)>0$ we have
\begin{equation}\label{eqn:dist-x-zi}
   d(x, z_i) \leq d(x, B_i) + \diam B_i + d(B_i, z_i) \leq (2+\alpha)d(B_i, Z).
 \end{equation}
Let $r$ be the smallest integer $\geq\frac{1+\alpha + \tau}{1-\tau}$ and, for each $i\in I$, define 
\begin{equation*}
 s_i:= \min\left\{s\in\Z: d(B_i, Z)\leq r^s\right\}.
\end{equation*}
For $x\in X\backslash Z$ and $i, j\in I$ with $\sigma_i(x)>0$ and $\sigma_j(x)>0$ we then have $|s_i-s_j|\leq 1$, as a consequence of \eqref{eqn:x-B-Z}.

Next, let $r$ be as above and set $\varepsilon':= 1$ and $\varepsilon:= [2\varepsilon' + 2(2+\alpha)]r$. Let $A, A'\subset J^k(\R^n)$ and $F, F'$ be as in \propref{prop:good-subset-Heis}. Choose nearest point projections $P: J^k(\R^n) \to A$ and $P': J^k(\R^n) \to A'$, necessarily discontinuous. For $i\in I$ let $P_i:= \delta_{r^{s_i}}\circ P \circ \delta_{r^{-s_i}}$ if $s_i$ is even, and $P_i:= \delta_{r^{s_i}}\circ P' \circ \delta_{r^{-s_i}}$ if $s_i$ is odd.  These maps are nearest-point projections to $\delta_{r^{s_i}}(A)$ and $\delta_{r^{s_i}}(A')$ respectively, so they satisfy
$$d_c(x,P_i(x))\le \varepsilon' r^{s_i}$$
for all $x$.
Define a map $\varphi: \Sigma^{(0)}\to J^k(\R^n)$ by
\begin{equation*}
 \varphi(e_i):=  P_i(f(z_i)).
\end{equation*}
We claim that for $x\in X\backslash Z$ and $i, j\in I$ with $\sigma_i(x)>0$ and $\sigma_j(x)>0$ we have
\begin{equation}\label{eqn:dist-h0}
 d_c(\varphi(e_i), \varphi(e_j)) \leq  [2\varepsilon' + 2(2+\alpha)] r^{\max\{s_i, s_j\}}.
\end{equation}
Indeed, 
\begin{align*}
  d_c(\varphi(e_i), \varphi(e_j)) & = d_c(P_i(f(z_i)), P_j(f(z_j)))\\
  &\le d_c(f(z_i), f(z_j))+\varepsilon' r^{s_i}+\varepsilon' r^{s_j}\\
  &\le d_c(z_i, z_j)+2 \varepsilon' r^{\max\{s_i, s_j\}}\\
  &\le d_c(z_i, x) +d_c(x, z_j)+2 \varepsilon' r^{\max\{s_i, s_j\}}\\
  &\le [2(2+\alpha)+2\varepsilon'] r^{\max\{s_i, s_j\}}
\end{align*}
where the last inequality is a consequence of \eqref{eqn:dist-x-zi} and the definition of $s_i$.  

\begin{figure*}[h]
 \begin{minipage}[b]{12cm}
 \psfrag{X}{$X$}
 \psfrag{Z}{$Z$}
  \psfrag{Si}{$\Sigma'$}
 \psfrag{J}{$J^k(\R^n)$}
 \psfrag{B1}{{\footnotesize $B_{i_1}$}}
 \psfrag{B2}{{\footnotesize $B_{i_2}$}}
 \psfrag{B3}{{\footnotesize $B_{i_3}$}}
 \psfrag{B4}{{\footnotesize $B_{i_4}$}}
  \psfrag{g}{$g$}
   \psfrag{1}{\tiny $1$}
 \psfrag{h}{\; $h$}
  \psfrag{e1}{\footnotesize $e_{i_1}$}
    \psfrag{e2}{\footnotesize $e_{i_2}$}
      \psfrag{e3}{\footnotesize $e_{i_3}$}
        \psfrag{e4}{\footnotesize $e_{i_4}$}
 \psfrag{fe1}{\footnotesize $\varphi(e_{i_1})$}
  \psfrag{fe2}{\footnotesize \; $\varphi(e_{i_2})$}
   \psfrag{fe3}{\footnotesize $\varphi(e_{i_3})$}
    \psfrag{fe4}{\footnotesize $\varphi(e_{i_4})$}
     \psfrag{x}{\tiny $x$}
  \psfrag{gx}{\tiny \;$g(x)$}
   \psfrag{fbx}{\footnotesize $\bar{f}(x)$}
    \psfrag{S1}{\footnotesize \; $S$} 
    \psfrag{S2}{\footnotesize \; $S'$}
     \psfrag{hS2}{\footnotesize \; $h(S')$}
      \psfrag{hS1}{\footnotesize \; $h(S)$}
      \psfrag{gB1}{\tiny \; $g(B_{i_1})$}
      \psfrag{gB2}{\tiny \; $g(B_{i_2})$}
         \psfrag{gB3}{\tiny \; $g(B_{i_3})$}
               \psfrag{gB4}{\tiny \; $g(B_{i_4})$}
        \psfrag{a}{\footnotesize  In the figure: $s_{i_1} = 0$, $s_{i_2}=s_{i_3}=1$, $s_{i_4}=2$}
         \psfrag{b}{ \footnotesize $\bullet$ $\varphi(e_{i_1}) =  P(f(z_{i_1}))\in A$}
          \psfrag{c}{ \footnotesize $\bullet$ $\varphi(e_{i_2}) = \delta_{r} \circ P' \circ \delta_{r^{-1}}(f(z_{i_2}))\in\delta_{r}(A')$}
           \psfrag{d}{ \footnotesize $\bullet$ $\varphi(e_{i_3}) = \delta_{r} \circ P' \circ \delta_{r^{-1}}(f(z_{i_3}))\in\delta_{r}(A')$}
           \psfrag{e}{ \footnotesize $\bullet$ $\varphi(e_{i_4}) = \delta_{r^2} \circ P \circ \delta_{r^{-2}}(f(z_{i_4}))\in\delta_{r^2}(A)$}
           \psfrag{E}{ \footnotesize $\bullet$ $h(S') = F([\varphi(e_{i_1}), \varphi(e_{i_2}), \varphi(e_{i_3})])$}
           \psfrag{F}{ \footnotesize $\bullet$ $h(S) = \delta_r(F'([\delta_{r^{-1}}(\varphi(e_{i_2})), \delta_{r^{-1}}(\varphi(e_{i_3})), \delta_{r^{-1}}(\varphi(e_{i_3}))]))$}
            \psfrag{G}{ \footnotesize $\bullet$ $h|_{S'} = F|$}
            \psfrag{H}{ \footnotesize $\bullet$ $h|_{S} = \delta_r\circ F'|$}
         \includegraphics*[width=12cm]{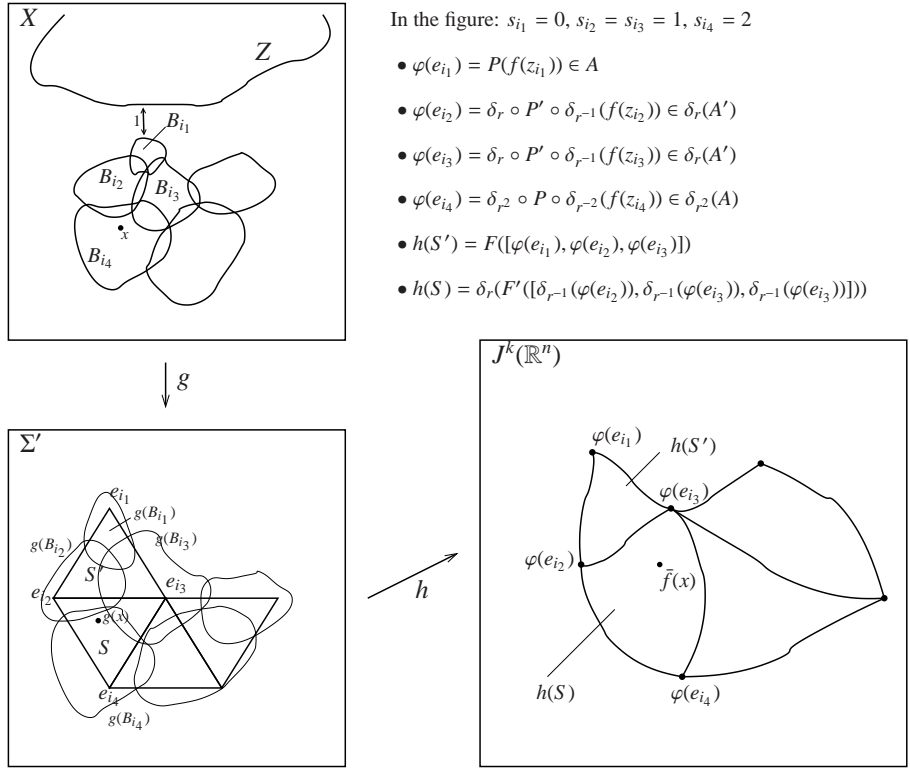}
         \caption{Construction of $g$ and $h$}\label{figure}
      \end{minipage}
\end{figure*}

We now define an extension $h:\Sigma'\to J^k(\R^n)$ of $\varphi$ as follows. Fix a simplex $S = [e_{i_1}, \dots, e_{i_m}]$ in $\Sigma'$. Set $t:= \min\{s_{i_1}, \dots, s_{i_m}\}$ and define $y_j:= \delta_{r^{-t}}(\varphi(e_{i_j}))$ for $j=1, \dots, m$. It follows from \eqref{eqn:dist-h0} that
\begin{equation*}
 d_c(y_{j_1}, y_{j_2})\leq [2\varepsilon' + 2(2+\alpha)]r=\varepsilon.
\end{equation*}
It is straightforward to check that $y_j\in\Gamma$ if $t$ is even and $y_j\in\Gamma'$ if $t$ is odd. Thus, $[e_{y_1}, \dots, e_{y_m}]$ is a simplex either in $\Sigma(\Gamma, \varepsilon, n)$ or in $\Sigma(\Gamma', \varepsilon, n)$, depending on whether $t$ is even or odd. We define $h$ such that
\begin{equation*}
 h|_S:= \left\{
  \begin{array}{ll}
  \delta_{r^t}\circ F|_{[e_{y_1}, \dots, e_{y_m}]} \circ \psi & \text{if $t$ is even}\\
 \delta_{r^t}\circ F'|_{[e_{y_1}, \dots, e_{y_m}]}  \circ \psi  & \text{if $t$ is odd,}
  \end{array}\right.
\end{equation*}
where $\psi$ is the natural isometry from $[e_{i_1}, \dots, e_{i_m}]$ onto $[e_{y_1}, \dots, e_{y_m}]$. The properties of $F$ and $F'$ guarantee that $h$ is well-defined and an extension of $\varphi$. We finally define an extension $\bar{f}: X\to J^k(\R^n)$ of $f$ so that, on $X\backslash Z$,
\begin{equation*}
 \bar{f} = h\circ g
\end{equation*}
and prove that $\bar{f}$ is Lipschitz continuous. For this, let first $x\in X\backslash Z$, and let $S$ be a simplex in $\Sigma'$ containing $g(x)$. Let $i\in I$ be such that $x\in B_i$, so $e_i$ is a vertex of $S$.  Now, if $z\in Z$ then we have
\begin{equation*}
 \begin{split}
   d_c(\bar{f}(x), f(z)) &\leq d_c(\bar{f}(x), \varphi(e_i)) + d_c(\varphi(e_i), f(z))\\
    &= d_c(h|_S(g(x)), h|_S(e_i)) + d_c(P_i(f(z_i)), f(z))\\
    &\leq C_1\varrho r^{s_i+1} + r^{s_i}\varepsilon' + d(z_i, z)\\
    &\leq (C_1\varrho r^2 + \varepsilon'r) d(B_i, Z) + d(z_i, x) + d(x, z)\\
    &<[C_1\varrho r^2 + \varepsilon'r + 2 + \alpha + 1]\, d(x, z),
 \end{split}
\end{equation*}
for a constant $C_1$ depending only on $n$. Note that $\varrho$ is the Lipschitz constant from \propref{prop:good-subset-Heis} and that we used \eqref{eqn:dist-x-zi} in the last inequality.
Now let $x, y\in X\backslash Z$ and let $S$ and $T$ be simplices in $\Sigma'$ containing $g(x)$ and $g(y)$, respectively. Pick $i,j\in I$ such that $x\in B_i$ and $y\in B_j$. Suppose first that $\sigma_i(y) = 0$ and $\sigma_j(x)=0$. Then $d(x,y)\geq d(y, B_i)\geq \tau d(B_i, Z) > \tau r^{s_i-1}$ and likewise $d(x, y)>\tau r^{s_j-1}$, so we obtain
\begin{equation*}
 \begin{split}
  d_c(\bar{f}(x), \bar{f}(y)) &\leq d_c(\bar{f}(x), \varphi(e_i)) + d_c(\bar{f}(y), \varphi(e_j)) + d_c(\varphi(e_i), \varphi(e_j))\\
  &\leq C_1\varrho [r^{s_i+1} + r^{s_j+1}] + 2\varepsilon' r^{\max\{s_i,s_j\}} + d(z_i,z_j)\\
  &\leq (2 r C_1\varrho  + 2\varepsilon')r^{\max\{s_i,s_j\}} + d(z_i,x) + d(x,y) + d(y,z_j)\\
  &\leq [2 r C_1\varrho  + 2\varepsilon' + 2(2+\alpha)]r^{\max\{s_i,s_j\}} + d(x,y)\\
  &\leq [2 r C_1\varrho  + 2\varepsilon' + 2(2+\alpha)]r \tau^{-1}d(x,y) + d(x,y)\\
 \end{split}
\end{equation*}
where the second-to-last inequality uses \eqref{eqn:dist-x-zi} and the last inequality uses the fact that $$d(x, y)>\tau r^{\max\{s_i,s_j\}-1}.$$

Finally, suppose that $\sigma_i(y)>0$ or $\sigma_j(x)>0$, for instance $\sigma_i(y)>0$. Then $e_i$ is a common vertex of $S$ and $T$. There thus exists a point $v\in S\cap T$ such that $d(g(x), v) + d(g(y), v) \leq \eta d(g(x), g(y))$ for some constant $\eta$ depending only on $n$. It follows that
\begin{equation*}
 \begin{split}
 d_c(\bar{f}(x), \bar{f}(y))&\leq d_c(h|_S(g(x)), h|_S(v)) + d_c(h|_T(g(y)), h|_T(v))\\
  &\leq \varrho r^{s_i+1} [d(g(x), v) + d(g(y),v)]\\
  &\leq \varrho \eta  r^2 d(B_i, Z)d(g(x), g(y)).
 \end{split}
\end{equation*}
There are at most $2n+1$ indices $k\in I$ with $\sigma_k(x)>0$ or $\sigma_k(y)>0$.
For all such $k\in I$ we compute, just as in \cite{Lang-Schlichenmaier},
\begin{equation*}
 \begin{split}
  \left|\frac{\sigma_k(x)}{\bar{\sigma}(x)} - \frac{\sigma_k(y)}{\bar{\sigma}(y)}\right| 
    &\leq  \left|\frac{\sigma_k(x)}{\bar{\sigma}(x)} - \frac{\sigma_k(y)}{\bar{\sigma}(x)}\right| 
       +  \left|\frac{\sigma_k(y)}{\bar{\sigma}(x)} - \frac{\sigma_k(y)}{\bar{\sigma}(y)}\right|\\
            &\leq \frac{1}{\bar{\sigma}(x)}\left(\left| \sigma_k(x) - \sigma_k(y)\right| + \left|\bar{\sigma}(x) - \bar{\sigma}(y)\right|\right)\\
    &\leq \frac{2(n+1)}{\bar{\sigma}(x)}d(x,y).
 \end{split}
\end{equation*}
It follows that $d(g(x), g(y))\leq\bar{c}d(x,y)/\bar{\sigma}(x)$ for some constant $\bar{c}$ only depending on $n$. Since $x\in B_i$, we have $\bar{\sigma}(x)\geq \sigma_i(x) = \tau d(B_i, Z)$, and we conclude that
\begin{equation*}
 d_c(\bar{f}(x), \bar{f}(y))\leq \varrho \eta  r^2 \bar{c} \tau^{-1}d(x,y).
\end{equation*}
This completes the proof.
\end{proof}

\end{document}